\documentclass[10pt]{article}

\usepackage{setspace}
\usepackage{float}
\usepackage{amsmath, amssymb}
\usepackage{ifthen}
\usepackage{mdframed} 
\usepackage{array}
\usepackage{listings}

\usepackage[inline]{enumitem}
\setlist{topsep=0pt, parsep=8pt}
\setlist[itemize]{topsep=0pt, itemsep=2pt, parsep=0pt}

\PassOptionsToPackage{hyphens}{url}
\PassOptionsToPackage{bookmarks,colorlinks,linkcolor=black,citecolor=black,pdfstartview=FitH,urlcolor=blue}{hyperref}
\usepackage{hyperref}
\hypersetup{pdfpagemode=UseNone}

\usepackage[rgb,table,usenames,dvipsnames]{xcolor}\convertcolorsUfalse
\usepackage{graphicx}

\usepackage{bookmark}

\definecolor{crimson}{HTML}{A41034}
\definecolor{crimsonlight}{HTML}{CC5500} 
\definecolor{crimsonmid}{HTML}{B8203D}   
\definecolor{crimsondark}{HTML}{780021}  
\definecolor{crimsongray}{HTML}{C4A4A9}  

\usepackage{array}
\usepackage{soul}
\usepackage[normalem]{ulem}

\usepackage{tikz}
\usetikzlibrary{
  matrix,
  decorations.pathmorphing,%
  decorations.markings,%
  decorations.pathreplacing,%
  decorations.text,%
  calc,%
  patterns,%
  intersections,
  backgrounds
}
\tikzset{>=latex} 
\tikzset{mark size=1.5pt, mark options=thin}
\newcommand\iscurrentenumi[1]{%
  \ifthenelse{\equal{#1}{\theenumi}}{}{#1}}
\setenumerate[1]{ref=\#\arabic*}
\setenumerate[2]{ref=\protect\iscurrentenumi{\theenumi}(\alph*)}
\newcommand\enumiiprefix[2]{%
  \ifthenelse{{\equal{#1}{\theenumi}}\AND{\equal{#2}{(\alph{enumii})}}}{}{}%
  \ifthenelse{{\equal{#1}{\theenumi}}\AND{\NOT{\equal{#2}{(\alph{enumii})}}}}{#2}{}%
  \ifthenelse{\equal{#1}{\theenumi}}{}{#1#2}%
}
\setenumerate[3]{ref=\protect\enumiiprefix{\theenumi}{(\alph{enumii})}\roman*}

\ifthenelse{\isundefined{\C}}{%
  \newcommand{\C}{\mathbb{C}}%
}{\renewcommand{\C}{\mathbb{C}}}
\ifthenelse{\isundefined{\R}}{%
  \newcommand{\R}{\mathbb{R}}%
}{\renewcommand{\R}{\mathbb{R}}}
\ifthenelse{\isundefined{\N}}{%
  \newcommand{\N}{\mathbb{N}}%
}{\renewcommand{\N}{\mathbb{N}}}

\ifthenelse{\isundefined{\F}}{%
  \newcommand{\F}{\mathbb{F}}%
}{\renewcommand{\F}{\mathbb{F}}}

\usepackage{amsthm}

\usepackage{amsthm}
\usepackage{xcolor}

\definecolor{codegray}{rgb}{0.9, 0.9, 0.9}
\definecolor{codeblue}{rgb}{0.0, 0.0, 0.6}

\lstdefinelanguage{Mathematica}{
  morekeywords={Module, Table, Do, If, Total, Length, Map, PartitionsP, IntegerPartitions, Tuples, Flatten, Function, Export, GraphicsGrid, MatrixPlot, Print},
  keywordstyle=\color{codeblue}\bfseries,
  sensitive=true,
  comment=[l]{(*},
  morecomment=[s]{(*}{*)},
  commentstyle=\color{gray}\ttfamily,
  stringstyle=\color{red},
}

\definecolor{theoremcolor}{rgb}{0.0, 0.2, 0.6} 
\definecolor{definitioncolor}{rgb}{0.0, 0.5, 0.0}
\definecolor{lemmacolor}{rgb}{0.5, 0.1, 0.1}     
\definecolor{propositioncolor}{rgb}{0.5, 0.3, 0.0}

\theoremstyle{definition}

\theoremstyle{plain}

\theoremstyle{remark}

\newmdenv[
  skipabove=0.5pt,
  skipbelow=10pt,
  leftmargin=0pt, 
  rightmargin=0pt,
  backgroundcolor=white!95,
  linecolor=black!20,
  linewidth=0.5pt,
  innertopmargin=1em,
  innerbottommargin=1em,
  innerleftmargin=0.6em,
  innerrightmargin=1em
]{proofsketchbox}

\newenvironment{proofsketch}
  {\begin{proofsketchbox}\textit{Proof.}}
  {\hfill$\square$\end{proofsketchbox}}

\numberwithin{equation}{section}
\numberwithin{theorem}{section}

\usepackage{tcolorbox} 

\newcommand{\ferrers}[1]{
    \begin{tikzpicture}[scale=0.5]
        \foreach \row[count=\i] in {#1} {
            \foreach \j in {1,...,\row} {
                \filldraw (\j,-\i) circle (0.2);
            }
        }
    \end{tikzpicture}
}

\newmdenv[
    linecolor=gray,         
    linewidth=2pt,          
    topline=false,          
    bottomline=false,       
    rightline=false,        
    skipabove=10pt,         
    skipbelow=10pt,         
    innerleftmargin=10pt,   
    innerrightmargin=0pt,   
    innertopmargin=5pt,     
    innerbottommargin=5pt,  
    backgroundcolor=white,  
]{example}

\title{Segre Characteristic Equivalence}
\begin{document}
\begin{center}
\textbf{\large SEGRE CHARACTERISTIC EQUIVALENCE}

\vspace{1ex}

\textit{Jessie Pitsillides} \\
\vspace{0.7ex}
\textit{\today}
\end{center}

\vspace{2ex}

Given only the dimension, $n$, of a square matrix $A \in M(n,\C)$, how many Segre Characteristic equivalent matrices are there? Jordan Normal Form Theorem states that any linear operator over $\C$ is similar to a matrix in Jordan Normal Form \cite[p.~324]{axler_linear_algebra}. 
As such, this is a question of counting the number of possible Jordan Normal Forms for a given dimension. So, equivalently, how many Jordan Normal Forms can an $n\times n$ matrix possibly have? 

\smallskip
\noindent
\textbf{1.\quad Definitions}

\noindent

We call two matrices $A, B$ in $M(n,\C)$ \textbf{similar} if there exists $S \in GL(n, \C)$ such that $B = S^{-1}AS$. A matrix's \textbf{Segre Characteristic} is the list of lists of integers that give the sizes of the blocks in the matrix's Jordan Canonical Form, grouped by eigenvalue. For example, the Segre Characteristic of
    \[
\begin{bmatrix}
\alpha & 1 & & & & & & & & \\
 &\alpha & & & & & & & &\\
 & &\alpha & & & & & & &\\
&  & &\beta & 1 & & & & &\\
&  & & &\beta & 1 & & & & \\
& & & & &\beta  & & & & \\
& & & & & &\gamma & & &\\
& & & & & & &\delta & 1 &\\
& & & &  & & & & \delta&\\
& & & &  & & & & &\delta
\end{bmatrix}
\]
is $[(2,1),3,1,(2,1)]$ \cite{FrazerDuncanCollar1963}. This notation groups block sizes associated with the same eigenvalue using parentheses, and wraps the entire list in square brackets. We call two matrices \textbf{Segre Equivalent} if they have the same Segre Characteristic.

Answering the counting question of how many Segre Equivalent matrices an $n\times n$ matrix may have is a question of partitions. A \textbf{partition} of a positive integer $n$ is a way of writing $n$ as a sum of positive integers, where the order of the summands does not matter. The \textbf{Partition Function}, $p(n)$, counts the number of partitions for $n$. The generating function of $p(n)$ is $\displaystyle\Pi_{k=1}^\infty(1-x^k)^{-1}$. The \textbf{Partition of Partition Function}, $P(n)$, counts how many partitions of partitions one can build. Its generating function is $\displaystyle\Pi_{k=1}^\infty(1-x^k)^{-p(k)}$. 

We also define \textbf{Rank Pattern}, \( R_{A}(\lambda) \), of \( A \) with respect to each \( \lambda \) to be the sequence $R_{A}(\lambda) := \{\text{rank}((A - \lambda I)^k) \mid k = 0, 1, 2, \dots, m\}$
          where \( m \) is the smallest integer such that \((A - \lambda I)^m = 0\), or where the rank stabilizes.

\smallskip
\noindent
\textbf{2.\quad Theorems}

\noindent
\textbf{Theorem 1.} \textit{The Rank Pattern explicitly determines the Segre Characteristic.}
\begin{proofsketch}
    Fix an eigenvalue, $\lambda_k$. Given that eigenvalue's Rank Pattern, $R_k$, we can calculate the nullity pattern, $Q_k$, which is given by $Q_k = n - R_k$. We then differentiate to obtain the nullity growth pattern, $q_k = Q_k - Q_{k-1}$, where $Q_0 = 0$. Finally, we compute the conjugate partition, $p$, of the nullity growth pattern. This partition gives the Segre Characteristic corresponding to that eigenvalue, in other words, the sizes of the Jordan Normal Form blocks corresponding to $\lambda_k$. 
\end{proofsketch}
The fact that the Rank Pattern of a matrix explicitly defines its Jordan structure is explicitly proven by Shafarevich and Remizov in 2012 \cite[p.~173]{Shafarevich}. 

Going from Rank Pattern to Segre Characteristic can be seen using \textbf{Ferrers Diagrams}; these are graphical representations of partitions \cite[p.~6]{Andrews1976p6}. The Segre Characteristic of a particular eigenvalue is given by the conjugate partition Ferrers Diagram of the nullity growth sequence corresponding to that eigenvalue. For example, suppose $A \in M(n,\C)$ has eigenvalue $\lambda$, and its rank pattern is $R_A(\lambda)=\{10,7,5,3,2,1,0\dots\}$. 
    So, its nullity growth sequence is: $\{3,2,2,1,1,1,0\dots\}$. This produces the Ferrers diagram: 
    \[
    \ferrers{3,2,2,1,1,1}
    \]
    Its conjugate transpose is:
    \[
    \ferrers{6,3,1}
    \]
    and thus the Segre Characteristic and hence Jordan block dimensions corresponding to $\lambda$ are $(6,3,1)$.

\textbf{Theorem 2.} \textit{There are $P(n)$ Segre Characteristic equivalent matrices of dimension $n$.}
\begin{proofsketch}
    The number of possible constituent Jordan blocks that an $n\times n$ matrix can be comprised of is the integer partitions of $n$. For each partition, we must count the number of ways we can distribute eigenvalues across the blocks. This is the number of unique subsets of the elements in the original partition. This is by definition the partitions of partitions of $n$, denoted $P(n)$.
\end{proofsketch}

For example, we can consider how many Segre Characteristic equivalent matrices there are for a $4 \times 4$ matrix. The number of possible Jordan block structures is the partitions of the integer $4$. Here are the $p(4) = 5$ possible structures: 
    \[
    \underbrace{\begin{bmatrix}
        &&&&\\
        &&&&\\
        &&&&\\
        &&&&
    \end{bmatrix}}_{[4]} 
    \underbrace{\begin{bmatrix}
        \begin{bmatrix}
            &&&\\
            &&&\\
            &&&
        \end{bmatrix}\\
    &\begin{bmatrix}
            &
        \end{bmatrix}
    \end{bmatrix}}_{[3,1]}
    \underbrace{\begin{bmatrix}
        \begin{bmatrix}
            &&\\
            &&
        \end{bmatrix}\\
        &\begin{bmatrix}
            &&\\
            &&
        \end{bmatrix}
    \end{bmatrix}}_{[2,2]}
    \underbrace{\begin{bmatrix}
        \begin{bmatrix}
            &&\\
            &&
        \end{bmatrix}\\
        &\begin{bmatrix}
            &
        \end{bmatrix}\\
        &&\begin{bmatrix}
            &
        \end{bmatrix}
    \end{bmatrix}}_{[2,1,1]}
    \underbrace{\begin{bmatrix}
        \begin{bmatrix}
            &
        \end{bmatrix}\\
        &\begin{bmatrix}
            &
        \end{bmatrix}\\
        &&\begin{bmatrix}
          &  
        \end{bmatrix}\\
        &&&\begin{bmatrix}
            &
        \end{bmatrix}
    \end{bmatrix}}_{[1,1,1,1]}
    \]

Now, we must consider all possibilities of how eigenvalues might be distributed across the Jordan blocks. The same eigenvalues can exist across multiple blocks, so we want to find the unique partitions of the sets of block dimensions above.

For $[4]$, we can only have one eigenvalue, so $[4]$ is the only possible Segre Characteristic. 
For $[3,1]$, we may have one eigenvalue, $[(3,1)]$, or two eigenvalues $[(3),(1)]$. For $[2,2]$, we may have one, $[(2,2)]$, or two $[(2),(2)]$ eigenvalues. For $[2,1,1]$, we may have one, $[(2,1,1)]$, two, $[(2,1), (1)]$ or $[(2), (1,1)]$, or three, $[(2),(1),(1)]$ eigenvalues. And finally for $[1,1,1,1]$, we may have one, $[(1,1,1,1)]$, two, $[(1,1,1),(1)]$, or $[(1,1),(1,1)]$, three, $[(1),(1),(1,1)]$, or four $[(1),(1),(1),(1)]$ eigenvalues. So overall, $P(4) = 1+2+2+4+5 = 14$, and so there are $14$ possible Segre Characteristic equivalent matrices for a $4 \times 4$ matrix.

We can obtain visualizations of the possible Jordan structures for $n \times n$ matrices.\footnote{Thank you to Oliver Knill for computer algebra assistance and for posing the original counting question, and the original rank pattern question with his colleague Poning Chen. Please see \hyperlink{AppA}{Appendix A} for Mathematica Code code to visualize all Segre possibilities given only $n$.} In these visualizations, the colored squares on the diagonals represent eigenvalues, and the black squares on superdiagonals are 1's. \\

\begin{figure}[H]
    \centering
    \includegraphics[width=0.6\textwidth]{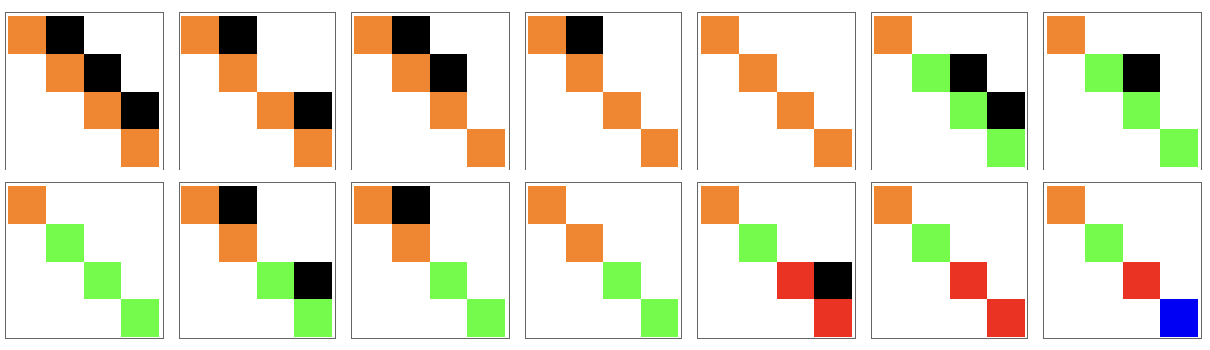}
    \caption{All 14 possible Segre Characteristics for $4 \times 4$ matrices}
    \label{fig:2segre4}
\end{figure}

\begin{figure}[H]
    \centering
    \includegraphics[width=0.6\textwidth]{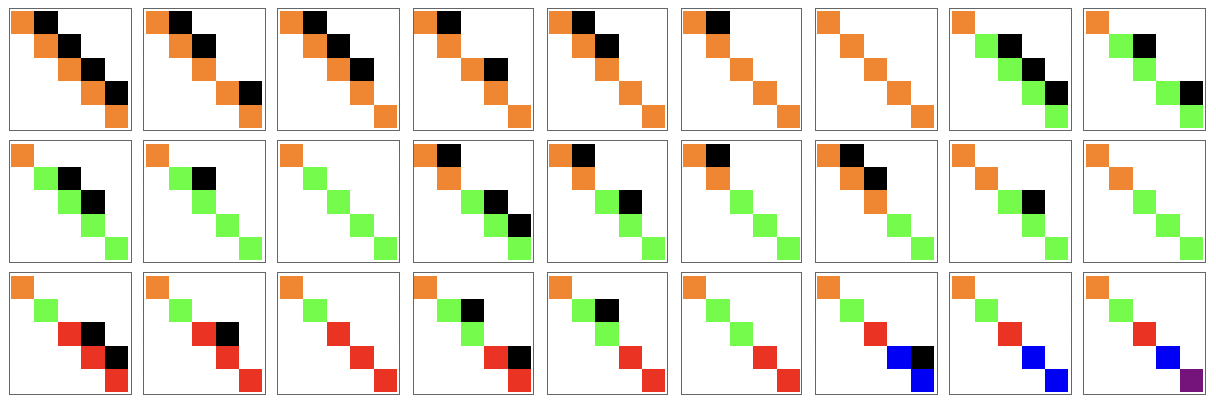} 
    \caption{All 27 possible Segre Characteristics for $5 \times 5$ matrices}
    \label{fig:2segre5}
\end{figure}

\begin{figure}[H]
    \centering
    \includegraphics[width=0.6\textwidth]{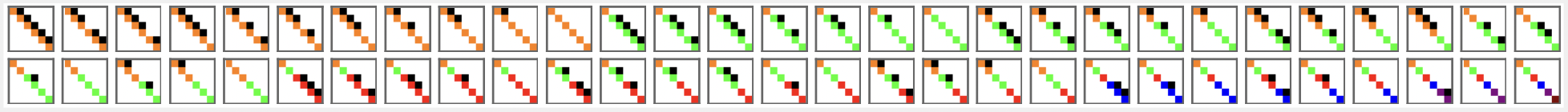} 
    \caption{All 58 possible Segre Characteristics for $6 \times 6$ matrices}
    \label{fig:2segre6}
\end{figure}

The sequence enumerating Segre Equivalence classes, known as \textit{partitions of partitions} (\href{https://oeis.org/search?q=A001970&language=english&go=Search}{OEIS Sequence A001970}) is
\[
1, 1, 3, 6, 14, 27, 58, 111, 223, 424, 817, 1527,\dots
\]
This sequence was originally published by J. J. Sylvester in 1854 in Philosophical Magazine \cite{Sylvester1854}, before Jordan Canonical form was formally developed in 1870. A historical account and tabulation of this sequence appears in Bromwich \cite[p.~60]{Bromwich}. 

\newpage
\smallskip
\noindent
\textbf{3.\quad Appendix}

\noindent
\hypertarget{AppA}{\textbf{A. \quad Mathematica Code for Segre Characteristic Visualization}}
\begin{lstlisting}
p[a_] := Product[PartitionsP[a[[k]]], {k, Length[a]}];
P[n_] := Total[Map[p, IntegerPartitions[n]]];
PP[m_] := 
  Module[{}, f[x_] = Product[1/(1 - x^k)^PartitionsP[k], {k, m + 1}]; 
   Last[CoefficientList[Series[f[x], {x, 0, m}], x]]];

BlockMatrix[A_] := 
  Module[{m = Length[A], a, n, B}, 
   a = Table[Length[A[[k]]], {k, Length[A]}];
   n = Total[a]; B = Table[0, {n}, {n}]; 
   s = Table[Sum[a[[j]], {j, k}], {k, 0, m}];
   Do[Do[
     B[[s[[l]] + i, s[[l]] + j]] = A[[l, i, j]], {i, a[[l]]}, {j, 
      a[[l]]}], {l, m}]; B];
M[a_, ev_] := 
  Module[{n = Total[a], m = Length[a], A, s}, 
   A = Table[If[i == j, ev, 0], {i, n}, {j, n}]; 
   s = Table[Sum[a[[j]], {j, k}], {k, 0, m}];
   Do[Do[A[[s[[l]] + k, s[[l]] + k + 1]] = 1, {k, 1, a[[l]] - 1}], {l,
      m}]; A];
M[s_] := BlockMatrix[Table[M[s[[k]], k + 1], {k, Length[s]}]];
Multipartitions[a_] := 
 Tuples[Table[IntegerPartitions[a[[k]]], {k, Length[a]}]]
Segre[n_] := Flatten[Map[Multipartitions, IntegerPartitions[n]], 1];
SortedSegre[n_] := Union[Map[Sort, Segre[n]]]; 
c = {White, Black, Orange, Green, Red, Blue, Purple};
Q[A_] := 
  MatrixPlot[A, ColorFunctionScaling -> False, Frame -> True, 
   FrameTicks -> False, 
   ColorFunction -> Function[x, c[[Floor[x] + 1]]]];

(* Using n=3 as an example input *)
S = SortedSegre[3]; J = Map[M, S];
GraphicsGrid[Partition[Map[Q, J], 3]]
\end{lstlisting}

\renewcommand{\refname}

\end{document}